# Free boundary minimal surfaces with connected boundary and arbitrary genus

Alessandro Carlotto, Giada Franz and Mario B. Schulz

**Abstract**

We employ min-max techniques to show that the unit ball in $\mathbb{R}^3$ contains embedded free boundary minimal surfaces with connected boundary and arbitrary genus.

## 1. Introduction

Over the last decade, the work by Fraser and Schoen [5–7] on extremals for Steklov eigenvalues has revitalised the study of free boundary minimal surfaces, whose origins go back at least to Courant. The theory has been developed in various interesting directions, yet many fundamental questions remain open. One of the most basic ones can be phrased as follows: does the unit ball of $\mathbb{R}^3$ contain free boundary minimal surfaces of any given genus $g \geq 0$ and any number of boundary components $b \geq 1$? In spite of significant advances, which we will survey below, the answer to such a question has proven to be quite elusive, even in special cases such as $g = b = 1$ (which is highlighted in *Open Question 1* of the recent survey [16]). There is a well-known analogy between the free boundary theory for the unit ball $B^3 \subset \mathbb{R}^3$ and the theory concerning closed minimal surfaces in the round three-dimensional sphere $S^3$: since for the latter Lawson proved in 1970 that there indeed exist in the sphere embedded minimal surfaces of arbitrary genus, there might be some reason to lean towards an affirmative answer. Here we fully solve the problem for the class of free boundary minimal surfaces with connected boundary:

**Theorem 1.1.** *For each $1 \leq g \in \mathbb{N}$ there exists an embedded free boundary minimal surface $M_g$ in $B^3$ with connected boundary, genus $g$ and dihedral symmetry $\mathbb{D}_{g+1}$.*

The dihedral group $\mathbb{D}_n$ is the symmetry group of a regular $n$-sided polygon. In our specific context, given $2 \leq n \in \mathbb{N}$ we shall define the *dihedral group* $\mathbb{D}_n$ of order $2n$ to be the subgroup of Euclidean isometries (acting on $\overline{B^3}$) generated by the rotation of angle $2\pi/n$ around the vertical axis $\xi_0 := \{(0,0,r) : r \in [-1,1]\}$ and by the rotations of angle $\pi$ around the $n$ horizontal axes $\xi_k := \{(r\cos(k\pi/n), r\sin(k\pi/n), 0) : r \in [-1,1]\}$ for $k \in \{1, \ldots, n\}$. We further define the *singular locus* of the dihedral group action to be $\mathcal{S} := \xi_0 \cup \xi_1 \cup \ldots \cup \xi_n$.

It follows from the proof of Theorem 1.1 that $M_g$ contains the horizontal axes $\xi_1, \ldots, \xi_{g+1}$. From a variational perspective, the surfaces in question are *unstable* (i.e. they have positive Morse index); in fact the main estimate in [1] implies that $M_g$ has index at least $\lceil 2g/3 \rceil$, thus growing (at least) linearly with the genus. In addition, these surfaces satisfy uniform, explicit lower and upper bounds on both their area and the length of their boundary curves (thanks to Corollary 3.8 below, and Theorem 2 in [21]). Although they exhibit some analogies with the higher-symmetry Chen–Gackstatter surfaces described in Section 5.5 of [11] (in particular: the same symmetry group), we note that, in a precise sense, the surfaces we construct *cannot* be regarded as the free boundary counterpart of known complete examples in $\mathbb{R}^3$, for indeed any complete embedded minimal surface with one end and finite total curvature must be a flat plane.

Differently from the approach presented by Lawson in [15] each surface $M_g$ is constructed by means of *global, variational methods*. More specifically, we employ the equivariant min-max theory developed by Ketover in [12] (for the closed case), and specified to the free boundary setting, with some striking applications, in [13]. In applying such machinery to prove Theorem 1.1 we first need, for any positive integer $g$, to carefully design a suitable genus $g$ equivariant sweepout so to ensure that not only the natural mountain-pass condition holds, but also (and more importantly) that the limit surface we obtain is attained with multiplicity *one*. This is a general issue that arises whenever one relies on min-max techniques, and it is in fact a rather delicate point in our construction. In turn, this aspect is crucial to make sure that both the number of boundary components and the genus are controlled throughout the process, i.e. as we take the limit of a min-max sequence.

We refer the reader to Section 3 of [16] for a broad overview of existence results, including those in higher-dimensional Euclidean balls or in the general setting of compact Riemannian manifolds with boundary, while we will focus here on the special case of $B^3$. There, the first non-trivial examples of (embedded) free boundary minimal surfaces, besides the flat disc and the critical catenoid, were obtained by Fraser and Schoen in [7]: these have genus zero and $b \geq 2$ boundary components. Through different methods, Folha–Pacard–Zolotareva constructed in [4] examples having genus zero or one and any *sufficiently large* number of boundary components. Later, Kapouleas and Li developed in [9] methods to desingularise the formal union of a disc and a critical catenoid to obtain free boundary minimal surfaces in $B^3$ with large genus and exactly three boundary components. Independently, as anticipated above, Ketover proposed a totally different approach to construct a sequence of surfaces that behaves, at least for large genus, exactly like the one in [9]. To get a pictorial description, this family can be regarded as a free boundary version of the Costa–Hoffman–Meeks minimal surfaces in $\mathbb{R}^3$. A different desingularisation scheme has been described, in the introduction of [10], to construct free boundary minimal surfaces having connected boundary and sufficiently large genus: such surfaces are obtained by regularising the intersection of two orthogonal discs in $B^3$ (by means of a suitable Scherk surface). Finally, in the same article Kapouleas and Wiygul constructed, essentially via perturbative methods, free boundary minimal surfaces in $B^3$ having connected boundary and prescribed high genus. Roughly speaking, what they presented is the base case for a more general procedure, that they call *stacking*, which consists in considering a certain number of parallel discs, joining them through suitable bridges and deforming the resulting objects in order to obtain novel free boundary minimal surfaces in the Euclidean ball.

So, to summarise, while interesting examples have been obtained in abundance, on the one hand the gluing/desingularisation methodologies do not (for their very nature) allow to obtain low-genus examples, and are only asymptotically effective, while on the other hand non-trivial technical obstacles arise if one aims at full topological control of min-max free boundary minimal surfaces. In this paper we focus on those issues and take care of them for the special sweepout family described in Section 2, which is enough to prove Theorem 1.1. The only general result we invoke, as an input for our main theorem, is the topological lower-semicontinuity theorem in [17] which implies, in our case, that the genus of $M_g$ is *at most $g$*. It turns out that, in order to conclude that equality holds for the genus, we first need to make sure that the boundary of $M_g$ is indeed connected. The question of controlling the number of boundary components in min-max constructions is notoriously delicate. In that respect, Li writes [17, pp. 324]: 'On the other hand, we note that it is impossible to get a similar bound on the connectivity (i.e., number of free boundary components) of the minimal surface.' In the context of the present paper, the conclusion that the boundary of $M_g$ must be connected is achieved through a rather surprising application of Simon's Lifting Lemma (cf. [3, Proposition 2.1]), which we present in Section 4. From there, the full control on the genus of the min-max surfaces we construct follows by the general characterisation of equivariant surfaces given in Appendix B.

The questions whether some of the surfaces we construct are unique for their given topological type, or whether they coincide (*for large genus*) with the family constructed in [10] via stacking methods remain open, as stands the related question whether $M_g$ could be characterised in terms of a maximising property for its first Steklov eigenvalue under the natural normalisation constraint.

**Acknowledgments.** The authors would like to express their gratitude to Ailana Fraser and Brian White for several enlightening conversations at the early stages of this project. This article was completed while A. C. was a visiting scholar at the *Institut Mittag-Leffler*: the excellent working conditions and the support of the *Royal Swedish Academy of Sciences* are gratefully acknowledged. The research of M. S. was funded by the EPSRC grant EP/S012907/1.

## 2. Effective sweepouts

**Definition 2.1** (cf. [13]). Let $\mathbb{D}_n$ be the dihedral group for some $n \geq 2$. A family $\{\Sigma_t\}_{t \in [0,1]}$ of closed subsets $\Sigma_t \subset \overline{B^3}$ with the following properties is called $\mathbb{D}_n$-sweepout of $B^3$.

(i) For all $t \in ]0,1[$ the set $\Sigma_t \subset \overline{B^3}$ is a smooth, embedded, compact surface with boundary $\partial \Sigma_t = \Sigma_t \cap \partial B^3$.

(ii) $\Sigma_0$ and $\Sigma_1$ are the union of a smooth, embedded, compact surface in $\overline{B^3}$ and a (possibly empty) finite collection of arcs in $\overline{B^3}$.

(iii) $\Sigma_t$ varies smoothly for $t \in ]0,1[$, and continuously, in the sense of varifolds, for $t \in [0,1]$.

(iv) Every $\Sigma_t$ is $\mathbb{D}_n$-equivariant, i.e. $\varphi(\Sigma_t) = \Sigma_t$ for all $\varphi \in \mathbb{D}_n$ and all $0 \leq t \leq 1$.

In this section we prove the following statement, one the existence of effective sweepouts of any genus.

**Lemma 2.2.** *Given $1 \leq g \in \mathbb{N}$ there exists a $\mathbb{D}_{g+1}$-sweepout $\{\Sigma_t\}_{t \in [0,1]}$ of $B^3$ such that $\mathscr{H}^2(\Sigma_0) = \mathscr{H}^2(\Sigma_1) = \pi$ and such that for every $0 < t < 1$*

- *the surface $\Sigma_t$ has genus $g$,*
- *the boundary of $\Sigma_t$ is connected,*
- *the area of $\Sigma_t$ is strictly less than $3\pi$.*

The idea behind our construction is to equivariantly glue *three* parallel discs through suitably controlled ribbons.

*Remark* 2.3. In a partly similar way Ketover [13, Theorem 5.1] glued two discs in order to variationally construct free boundary minimal surfaces of genus zero and $b \geq 2$ boundary components (to be compared with the existence result by Fraser and Schoen [7, Theorem 1.1]).

Let $1 \leq g \in \mathbb{N}$ be fixed. Let $D = \{x \in \overline{B^3} : x_3 = 0\}$ be the equatorial disc in the closed unit ball and let $\overline{B_\varepsilon}(p) = \{x \in \mathbb{R}^3 : |x-p| \leq \varepsilon\}$ denote the closed ball of radius $\varepsilon > 0$ around any given $p \in \mathbb{R}^3$. For all $k \in \{0, \ldots, g\}$ we consider the points

$$p_k^\pm := \left(\cos\left(\frac{2k \pm \frac{1}{2}}{g+1}\pi\right), \sin\left(\frac{2k \pm \frac{1}{2}}{g+1}\pi\right), 0\right) \tag{1}$$

on the equator and the subsets

$$D_\varepsilon^\pm := D \setminus \bigcup_{k=0}^{g} B_\varepsilon(p_k^\pm), \tag{2}$$

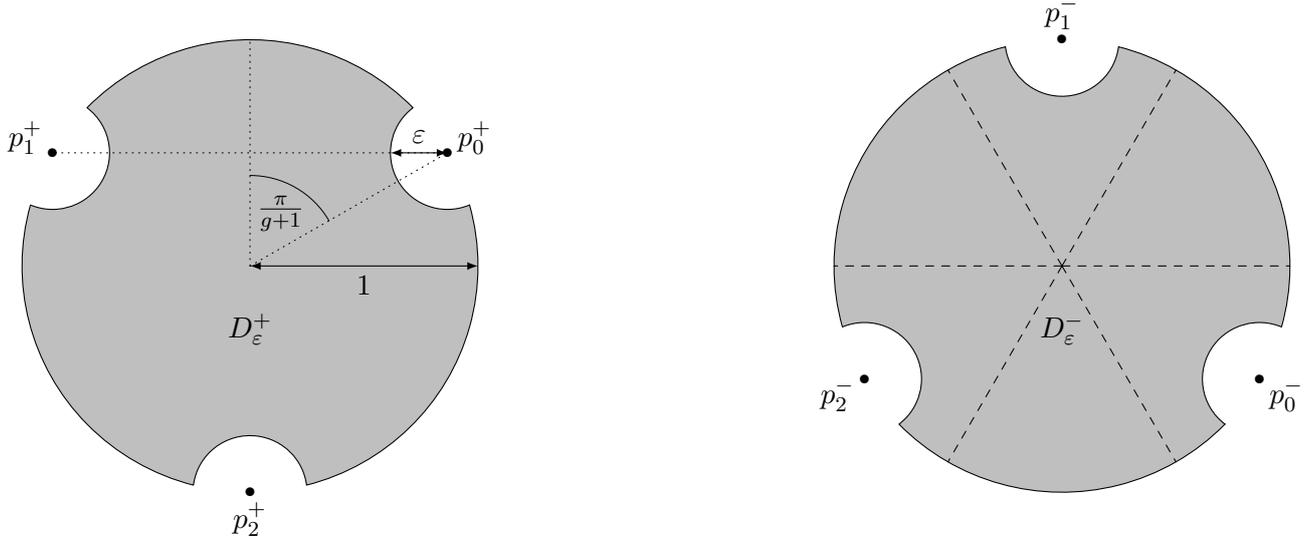

Figure 1: Top view of the sets $D_\varepsilon^\pm$ for $\varepsilon = 1/4$ and $g = 2$.

as shown in Fig. 1, which we then scale and translate upwards (or downwards) to define

$$D_{t,\varepsilon}^\pm := \bigl(\sqrt{1-t^2} D_\varepsilon^\pm\bigr) \pm (0,0,t) \tag{3}$$

for all $t \in [0,1]$. Now we connect the three sets $D_{t,\varepsilon}^+$, $D_{t,\varepsilon}^-$ and $(D_\varepsilon^+ \cap D_\varepsilon^-)$ in a $\mathbb{D}_{g+1}$-equivariant way by means of $2(g+1)$ ribbons. Note that each of these sets is a non-empty, connected subset of $\overline{B^3}$ provided that $\varepsilon < \sin(\pi/(2g+2))$. Let $0 < t_0 < 1$ be a fixed, small value which will be specified later in (11). For each $t \in [t_0, 1]$ we define

$$\Omega_{t,\varepsilon}^\pm := \bigcup_{s \in [0,t]} D_{s,\varepsilon}^\pm, \qquad S_{t,\varepsilon}^\pm := \overline{\partial \Omega_{t,\varepsilon}^\pm \setminus (\partial B^3 \cup D)}, \qquad \Sigma_t := S_{t,\varepsilon}^+ \cup S_{t,\varepsilon}^- \cup (D_\varepsilon^+ \cap D_\varepsilon^-). \tag{4}$$

In (4) the symbol $\partial$ refers to the topological boundary in $\mathbb{R}^3$. Moreover, we allow $\varepsilon \colon [t_0, 1[ \to ]0, \varepsilon_0]$ to be a continuous function of $t$, bounded from above by some sufficiently small $\varepsilon_0 > 0$ which we choose later in (11) depending on $t_0$ and $g$, such that $\varepsilon(t) \to 0$ as $t \nearrow 1$. Then we define $\Sigma_1$ to be the union of the equatorial disc $D$ with the (shortest) geodesic arcs connecting $p_k^+$ with the north pole and $p_k^-$ with the south pole for each $k \in \{0, \dots, g\}$. The construction is visualised in the first, second and third image of Fig. 2.

Arriving at $\Sigma_{t_0}$, one would like to increase $\varepsilon$ (as we further *decrease* $t$) in order to retract the three sets $D_{t,\varepsilon}^+$, $D_{t,\varepsilon}^-$ and $(D_\varepsilon^+ \cap D_\varepsilon^-)$ as illustrated in the fourth image of Fig. 2. However, this requires refined control on the area of the widening ribbons as we point out in the following statement.

**Lemma 2.4.** *Let $\Sigma_t$ be as given in* (4). *Then its area satisfies*

$$\mathscr{H}^2(\Sigma_t) \leq 3\pi - 2\pi\bigl(t^2 - (g+1)\varepsilon t\bigr). \tag{5}$$

*Proof.* Obviously, the set $S_{t,\varepsilon}^+$ has the same area as $S_{t,\varepsilon}^-$, by symmetry. Furthermore $S_{t,\varepsilon}^+$ is the union of $D_{t,\varepsilon}^+$ defined in (3) with $(g+1)$ ribbons. By construction,

$$\mathscr{H}^2(D_{t,\varepsilon}^+) \leq \pi(1-t^2), \qquad\qquad \mathscr{H}^2(D_\varepsilon^+ \cap D_\varepsilon^-) \leq \pi. \tag{6}$$

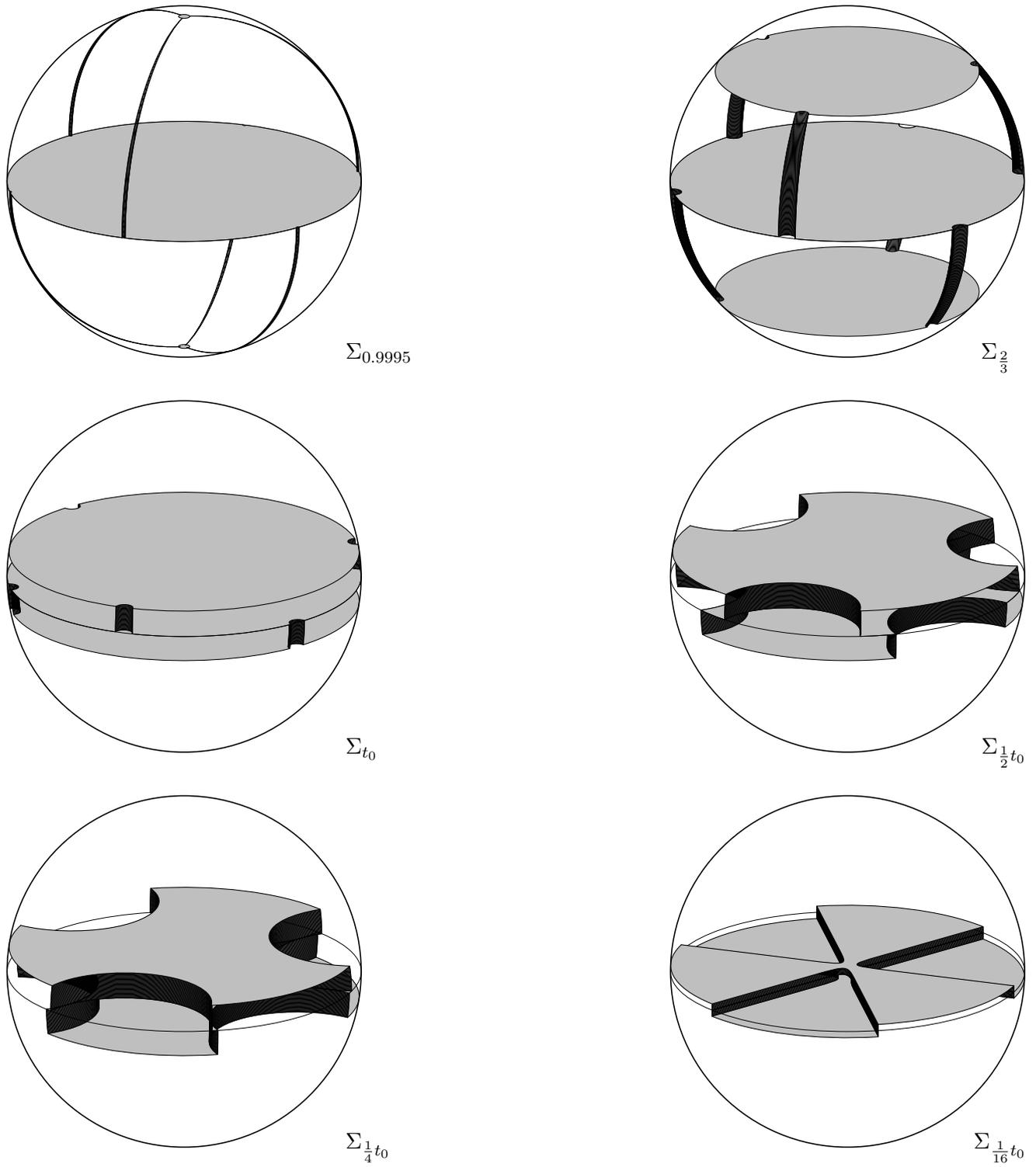

Figure 2: Construction of an effective sweepout in the case $g = 2$. In the first three images, $\varepsilon$ has been increased and relation (11) ignored for the sake of clarity.

The intersection of one ribbon with the horizontal plane at height $s \in [0,t]$ is an arc of length less than $\pi\sqrt{1-s^2}\varepsilon$. Hence, using the coarea formula one gets at once that the area of one ribbon is bounded from above by

$$\int_0^t \sqrt{1 + \frac{s^2}{1-s^2}} \sqrt{1-s^2}\, \pi\varepsilon\, ds = \pi\varepsilon t.$$

Therefore, $\mathscr{H}^2(\Sigma_t) = 2\mathscr{H}^2(S_{t,\varepsilon}^+) + \mathscr{H}^2(D_\varepsilon^+ \cap D_\varepsilon^-) \leq 2\pi(1-t^2) + 2(g+1)\pi\varepsilon t + \pi$, which allows to conclude by simply rearranging the terms. □

*Remark* 2.5. Lemma 2.4 implies that if $\varepsilon > 0$ is *small* compared to $t_0$, more precisely, if $(g+1)\varepsilon < t_0$, then $\mathscr{H}^2(\Sigma_t) < 3\pi$ holds for all $t \in [t_0, 1]$. In estimate (6) we did not take into account that by definition (2) of $D_\varepsilon^\pm$, small balls of radius $\varepsilon > 0$ are removed around the points $p_k^\pm$ (and thus similarly for $D_{t,\varepsilon}^\pm$). If we subtracted these contributions, one could then easily prove that the inequality $\mathscr{H}^2(\Sigma_{t_0}) < 3\pi$ would also hold for $\varepsilon > 0$ *large* compared to $t_0$. However, even the improved right-hand side of (5) would *not* stay below $3\pi$ if we increased $\varepsilon$ continuously from small to large values. For instance, that bound would be violated if one took $\varepsilon = 2t_0/(g+1)$. This is the reason why we need to refine the construction for $t < t_0$ and appeal to the so-called *catenoid estimate* instead.

At this stage, it would be possible to proceed by appealing to a suitable variant (for boundary points) of [14, Theorem 2.4]. However, for our specific scopes we will work out the explicit construction in our Euclidean setting.

Fix $0 < r < \sin(\pi/(2g+2))$ and $0 < h < \min\{\tanh(1)/2, 1/5\}r = r/5$. Moreover we choose $h$ such that we also have that $-\log h > 8(g+1)$. For every $s \geq 0$ consider the surfaces

$$C_s^{r,h} := \left\{ x \in \mathbb{R}^3 \ : \ \sqrt{x_1^2 + x_2^2} = \frac{r\cosh(sx_3)}{\cosh(sh)},\ |x_3| \leq h \right\}, \tag{7}$$

which all span two parallel circles of radius $r$ and distance $2h$. If $s$ is chosen such that

$$rs = \cosh(sh), \tag{8}$$

then $C_s^{r,h}$ is a subset of a (rescaled) catenoid and hence a minimal surface. For our choices of $r$ and $h$, equation (8) has two positive solutions $s_1(r,h) < s_2(r,h)$. The smaller one corresponds to the stable catenoid and the larger one to the unstable catenoid.

The family of surfaces in question interpolates between the cylinder at $s = 0$ and the union of a line segment with two discs of radius $r$ which we denote by $C_\infty^{r,h}$. The unstable catenoid can be regarded as the slice of largest area in this family, as we prove in Appendix A. The catenoid estimate given in [14, Proposition 2.1], combined with Lemma A.1, implies that we can choose $h$ possibly smaller (only depending on $r$) such that for all $s \geq 0$

$$\mathscr{H}^2(C_s^{r,h}) \leq \mathscr{H}^2(C_\infty^{r,h}) + \frac{4\pi h^2}{(-\log h)} = 2\pi r^2 + \frac{4\pi h^2}{(-\log h)}. \tag{9}$$

Let $E_R := \{x \in \mathbb{R}^3 \ : \ x_1^2 + (x_2 - R)^2 + x_3^2 < R^2\}$ be the ball of radius $R > 1$ around the point $(0, R, 0)$. By symmetry, the catenoid estimate holds under restriction to the half-space $E_\infty := \{x \in \mathbb{R}^3 \ : \ x_2 > 0\}$, i.e.

$$\sup_{s \geq 0} \mathscr{H}^2(C_s^{r,h} \cap E_\infty) \leq \mathscr{H}^2(C_\infty^{r,h} \cap E_\infty) + \frac{2\pi h^2}{(-\log h)}.$$

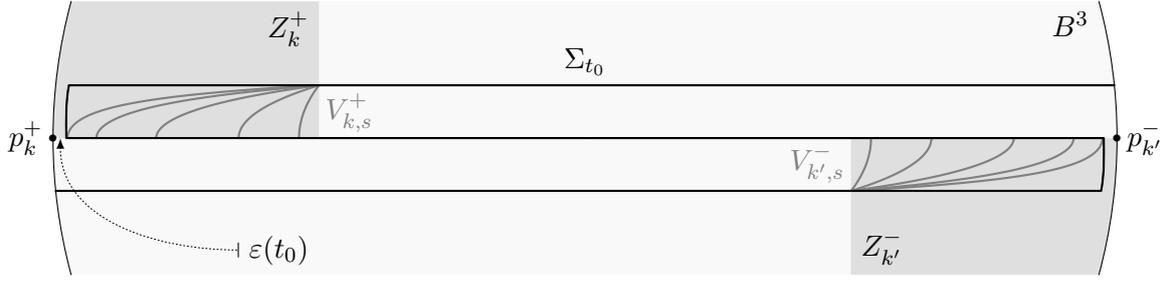

Figure 3: Replacement of $\Sigma_{t_0} \cap Z_k^\pm$ with $V_{k,s}^\pm$.

By a simple continuity argument, there exists $R_0 = R_0(r,h) > 1$ such that

$$\sup_{s \geq 0} \mathscr{H}^2(C_s^{r,h} \cap E_{R_0}) \leq \mathscr{H}^2(C_\infty^{r,h} \cap E_{R_0}) + \frac{4\pi h^2}{(-\log h)}.$$

Hence, renaming $r \mapsto r/R_0$ and $h \mapsto h/R_0$ (which corresponds to rescaling the whole picture by a factor of $1/R_0$), we obtain

$$\sup_{s \geq 0} \mathscr{H}^2(C_s^{r,h} \cap E_1) \leq \mathscr{H}^2(C_\infty^{r,h} \cap E_1) + \frac{4\pi h^2}{(-\log h + \log R_0)} \leq \mathscr{H}^2(C_\infty^{r,h} \cap E_1) + \frac{4\pi h^2}{(-\log h)}. \tag{10}$$

Observe that the conditions imposed on the smallness of $r$, $h$ and $h/r$ are still fulfilled.

We recall that $\Sigma_t$ was defined for all $t \in [t_0, 1]$ in (4), where we are free to choose first $t_0 > 0$ and then $\varepsilon_0 > 0$ such that

$$t_0 = h, \qquad \varepsilon_0 = \frac{t_0}{2(g+1)}. \tag{11}$$

By Lemma 2.4, this choice for $\varepsilon_0$ ensures that for all $t \in [t_0, 1]$

$$\mathscr{H}^2(\Sigma_t) \leq (3 - t_0^2)\pi. \tag{12}$$

For each $k \in \{0, \ldots, g\}$, let

$$Z_k^+ := \{x \in \overline{B^3} : \operatorname{dist}((x_1, x_2, 0), p_k^+) < r, \ x_3 > 0\},$$
$$Z_k^- := \{x \in \overline{B^3} : \operatorname{dist}((x_1, x_2, 0), p_k^-) < r, \ x_3 < 0\}$$

as shown in Fig. 3 and in Fig. 4 on the left. Note that, by the very way we have defined our parameters it follows that $r > 5t_0 \geq 10(g+1)\varepsilon_0$. We shall now replace $\Sigma_{t_0} \cap Z_k^+$ with a copy of the upper half of the surface $C_s^{r,t_0} \cap E_1$ after a suitable horizontal translation and rotation mapping $0 \mapsto p_k^+$ and $E_1 \mapsto B^3$. Similarly, $\Sigma_{t_0} \cap Z_k^-$ is replaced by a copy of the lower half of $C_s^{r,t_0} \cap E_1$. We denote those copies by $V_{k,s}^+$ and $V_{k,s}^-$, respectively. Initially, we choose $s = s_0$ such that

$$\frac{r}{\cosh(s_0 t_0)} = \varepsilon(t_0)$$

which ensures a continuous gluing of $V_{k,s_0}^\pm$ and $D_\varepsilon^+ \cap D_\varepsilon^- \subset \Sigma_{t_0}$ at height $x_3 = 0$. Moreover, assuming that $\varepsilon(t_0) \in {]0, \varepsilon_0[}$ is sufficiently small (so that $s_0$ will be very large), the surfaces $\Sigma_{t_0} \cap Z_k^\pm$ and $V_{k,s_0}^\pm$ are arbitrarily close such that we can continuously deform $\Sigma_{t_0} \cap Z_k^+$ into $V_{k,s_0}^+$ without significantly increasing the area. Then, as $t$ decreases further from $t_0$ to $t_0/2$, we decrease $s$ from $s_0$ to 0 and define $\Sigma_t$ through similar gluings of $V_{k,s}^\pm$ and $\Sigma_{t_0} \setminus Z_k^\pm$ as shown in Fig. 3 and Fig. 2, third and forth image. By (10) and by (12), we have for all $t \in [t_0/2, t_0[$

$$\mathscr{H}^2(\Sigma_t) \leq \mathscr{H}^2(\Sigma_{t_0}) - (g+1)\mathscr{H}^2(C_\infty^{r,t_0} \cap E_1) + (g+1)\mathscr{H}^2(C_s^{r,t_0} \cap E_1)$$
$$\leq (3 - t_0^2)\pi + (g+1)\frac{4\pi t_0^2}{(-\log t_0)} < 3\pi, \tag{13}$$

the last inequality relying on the fact that $-\log t_0 > 8(g+1)$. Now, observe further that

$$\mathscr{H}^2(\Sigma_{t_0/2}) \leq 3\pi - 2(g+1)\left(\frac{\pi}{2}r^2 - 2\pi r t_0\right) = 3\pi - (g+1)\pi(r^2 - 4rt_0) < 3\pi, \tag{14}$$

where we have used that $r > 5h = 5t_0$. It is now easy to see that it is possible to define $\Sigma_t$ for $t \in [0, t_0/2]$ in such a way $\mathscr{H}^2(\Sigma_t)$ is decreasing as $t$ decreases and $\Sigma_0$ is the equatorial disc. Indeed, thanks to (14) we see at once that by increasing $r$ till the threshold value $\bar{r} = \sin(\pi/(2g+2))$, so by removing larger discs as we vary $t \in [t_0/4, t_0/2]$ gives an area-decreasing deformation; then, for $t \in [0, t_0/4]$ we can just perform a simple retraction (see Fig. 4 and Fig. 2, fifth and sixth image).

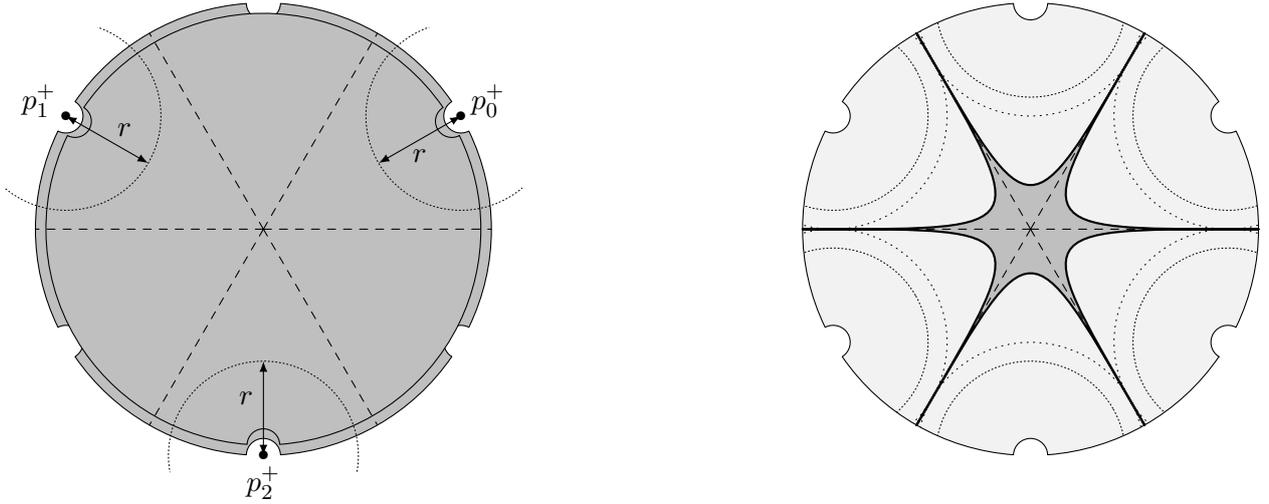

Figure 4: Implementing the catenoid estimate and making a further retraction.

*Proof of Lemma 2.2.* Let $\{\Sigma_t\}_{t \in [0,1]}$ be as constructed above. We define the desired smooth surfaces by regularising $\Sigma_t$ for all $0 < t < 1$ equivariantly (without renaming), a process which can be performed without violating the strict $3\pi$ upper bound on the area. We note that at $t = 0$, the origin is a singular point, where the genus of $\Sigma_t$ collapses as $t \searrow 0$. However, for $0 < t < 1$, we obtain a smooth family of genus $g$ surfaces, as claimed. □

## 3. Saturation of the sweepout and its width

In order to apply a min-max procedure, we need to consider a saturation of the sweepout given by Lemma 2.2, as in [13, Section 3].

**Definition 3.1.** A smooth map $\Phi \colon [0,1] \times \overline{B^3} \to \overline{B^3}$ is said to be a $\mathbb{D}_n$-*isotopy* for some $2 \leq n \in \mathbb{N}$ if

(i) $\Phi_t := \Phi(t, \cdot)$ is a diffeomorphism of $\overline{B^3}$ for all $0 \leq t \leq 1$;

(ii) $\Phi_0$ and $\Phi_1$ coincide with the identity map in $B^3$;

(iii) $\varphi \circ \Phi_t = \Phi_t \circ \varphi$ for all $0 \leq t \leq 1$ and all $\varphi \in \mathbb{D}_n$.

**Definition 3.2.** Given the $\mathbb{D}_{g+1}$-sweepout $\{\Sigma_t\}_{t\in[0,1]}$ constructed in Lemma 2.2, we define its $\mathbb{D}_{g+1}$-*saturation* as

$$\Pi := \{\{\Phi_t(\Sigma_t)\}_{t\in[0,1]} \ : \ \Phi \colon [0,1] \times \overline{B^3} \to \overline{B^3} \text{ is a } \mathbb{D}_{g+1}\text{-isotopy}\}.$$

Then the *min-max width* of $\Pi$ is defined as

$$W_\Pi := \inf_{\{\Lambda_t\}\in\Pi} \sup_{t\in[0,1]} \mathscr{H}^2(\Lambda_t).$$

In this section we prove that, in our context, the min-max width is strictly larger than $\pi$. For this purpose, it is helpful to introduce some terminology about finite perimeter sets (the reader is referred to e.g. [18, Chapter 12]).

Hereafter, any subset of the form $\{x \in B^3 \ : \ x \cdot v \geq 0\}$ for some $v \in \mathbb{R}^3$ is called a *half-ball*.

**Definition 3.3.** We say that a finite perimeter set $E \subset B^3$ is $\mathbb{D}_n$-equivariant if, for all $\varphi \in \mathbb{D}_n$, the set $\varphi(E)$ coincides either with $E$ or $B^3 \setminus E$ up to a negligible set.

**Lemma 3.4.** *Let $\{\Sigma_t\}_{t\in[0,1]}$ be the $\mathbb{D}_{g+1}$-sweepout in Lemma 2.2 and let $\Pi$ be its $\mathbb{D}_{g+1}$-saturation. Then for every $\{\Lambda_t\}_{t\in[0,1]} \in \Pi$ there exists a family $\{F_t\}_{t\in[0,1]}$ of $\mathbb{D}_{g+1}$-equivariant finite perimeter sets such that the following properties hold.*

(i) *$F_0$ is the upper half-ball and $F_1$ is the lower half-ball.*

(ii) *The family $\{F_t\}_{t\in[0,1]}$ is continuous in the sense of finite perimeter sets, i.e. $\mathscr{H}^3(F_t \triangle F_{t_0}) \to 0$ whenever $t \to t_0$, where $F_t \triangle F_{t_0} := (F_t \setminus F_{t_0}) \cup (F_{t_0} \setminus F_t)$.*

(iii) *The finite perimeter sets $F_t$ are $\mathbb{D}_{g+1}$-equivariant for all $0 \leq t \leq 1$.*

(iv) *For every $0 \leq t \leq 1$, $\Lambda_t$ is the relative boundary of $F_t$ in $B^3$; namely $\Lambda_t \setminus \partial B^3 = \partial F_t \setminus \partial B^3$.*

(v) *The volume of $F_t$ is half the volume of $B^3$ for all $0 \leq t \leq 1$, i.e. $\mathscr{H}^3(F_t) = \mathscr{H}^3(B^3)/2$ for all $0 \leq t \leq 1$.*

*Proof.* By the construction of the sweepout $\{\Sigma_t\}_{t\in[0,1]}$ in Lemma 2.2, we easily obtain that there exists a family $\{F_t^\Sigma\}_{t\in[0,1]}$ of $\mathbb{D}_{g+1}$-equivariant finite perimeter sets with properties (i)–(v). In particular, we can choose $\{F_t^\Sigma\}_{t\in[0,1]}$ such that $F_0^\Sigma$ is the upper half-ball and $F_1^\Sigma$ is the lower half-ball.

Now let us consider any other sweepout $\{\Lambda_t\}_{t\in[0,1]} \in \Pi$. By definition of saturation there exists a $\mathbb{D}_{g+1}$-isotopy $\Phi \colon [0,1] \times \overline{B^3} \to \overline{B^3}$ such that $\Lambda_t = \Phi_t(\Sigma_t)$ for all $0 \leq t \leq 1$. We want to prove that $\{F_t\}_{t\in[0,1]}$ defined by $F_t := \Phi_t(F_t^\Sigma)$ is a family of $\mathbb{D}_{g+1}$-equivariant finite perimeter sets as required in the statement. Observe that properties (i), (ii), (iii) and (iv) are straightforward, hence we just need to check (v).

Let us consider $\psi \in \mathbb{D}_{g+1}$ given by the rotation (of angle $\pi$) around the isotropy axis $\xi_1$ and observe that $\psi(F_t^\Sigma) = B^3 \setminus F_t^\Sigma$ for every $t \in [0,1]$. Hence we have that

$$\psi(F_t) = \psi(\Phi_t(F_t^\Sigma)) = \Phi_t(\psi(F_t^\Sigma)) = \Phi_t(B^3 \setminus F_t^\Sigma) = B^3 \setminus \Phi_t(F_t^\Sigma) = B^3 \setminus F_t,$$

which proves $\mathscr{H}^3(F_t) = \mathscr{H}^3(B^3)/2$ and concludes the proof. $\square$

We now denote by $\Xi$ the set of all families of $\mathbb{D}_{g+1}$-equivariant finite perimeter sets for which properties (i)–(v) in Lemma 3.4 hold for some $\mathbb{D}_{g+1}$-sweepout $\{\Lambda_t\}_{t\in[0,1]} \in \Pi$. Since there is a one-to-one correspondence between $\Pi$ and $\Xi$, we derive the following conclusion.

**Corollary 3.5.** *Under the hypotheses of the previous lemma, we have that*

$$W_\Pi = \inf_{\{F_t\}\in\Xi} \sup_{t\in[0,1]} P(F_t; B^3),$$

*where $P(F_t; B^3)$ denotes the relative perimeter of the finite perimeter set $F_t$ in $B^3$.*

In order to prove that $W_\Pi$ is in fact strictly larger than $\pi$, we first need the following stability lemma for the isoperimetric inequality.

**Lemma 3.6** (Stability of the isoperimetric inequality). *Fix $2 \leq n \in \mathbb{N}$. Then, for every $\varepsilon > 0$ there exists $\delta > 0$ such that, given a $\mathbb{D}_n$-equivariant finite perimeter set $F \subset B^3$ with Lebesgue measure $\mathscr{H}^3(F) = \mathscr{H}^3(B^3)/2$ and relative perimeter $P(F; B^3) \leq \pi + \delta$, there exists a $\mathbb{D}_n$-equivariant half-ball $\tilde{F}$ with $\mathscr{H}^3(F \triangle \tilde{F}) \leq \varepsilon$.*

*Proof.* Towards a contradiction, assume that there exist $\varepsilon > 0$ and a sequence $\{F_k\}_{k\in\mathbb{N}}$ of $\mathbb{D}_n$-equivariant finite perimeter sets satisfying $\mathscr{H}^3(F_k) = \mathscr{H}^3(B^3)/2$ and $P(F_k; B^3) \leq \pi + \delta_k$ for $\delta_k \to 0$ as well as $\mathscr{H}^3(F_k \triangle \tilde{F}) \geq \varepsilon$ for every $\mathbb{D}_n$-equivariant half-ball $\tilde{F}$.

By the compactness theorem for finite perimeter sets (see [18, Theorem 12.26]), there exists a $\mathbb{D}_n$-equivariant finite perimeter set $F_\infty \subset B^3$ such that a subsequence of $\{F_k\}_{k\in\mathbb{N}}$, which we do not rename, satisfies $\mathscr{H}^3(F_\infty \triangle F_k) \to 0$ as $k \to \infty$. In particular we have that $\mathscr{H}^3(F_\infty) = \mathscr{H}^3(B^3)/2$. Moreover, by lower semicontinuity of the perimeter ([18, Proposition 12.15]), it holds

$$P(F_\infty; B^3) \leq \liminf_{k\to\infty} P(F_k; B^3) = \pi.$$

Hence, $F_\infty$ is a finite perimeter set in $B^3$ with $\mathscr{H}^3(F_\infty) = \mathscr{H}^3(B^3)/2$ and $P(F_\infty; B^3) \leq \pi$, which implies that $F_\infty$ is a half-ball. Indeed, the reduced boundary $\partial^* F_\infty \cap B^3$ of $F_\infty$ in $B^3$ is smooth analytic with constant mean curvature by Theorem 27.4 in [18] (see pp. 386–389 therein for historical notes) and thus it is an equatorial disc by [2, Satz 1] (see also [20, Theorem 5]). However, this contradicts the choice of the sequence $\{F_k\}_{k\in\mathbb{N}}$ and concludes the proof. □

**Proposition 3.7.** *Fix $2 \leq n \in \mathbb{N}$. Then, there exists $\delta_0 > 0$ with the following property. Let $\{F_t\}_{t\in[0,1]}$ be a family of $\mathbb{D}_n$-equivariant finite perimeter sets in the unit ball $B^3$ such that*

(i) *$\{F_t\}_{t\in[0,1]}$ is continuous in the sense of finite perimeter sets, i.e. $\mathscr{H}^3(F_t \triangle F_{t_0}) \to 0$ whenever $t \to t_0$;*

(ii) *$\mathscr{H}^3(F_t) = \mathscr{H}^3(B^3)/2$ for all $0 \leq t \leq 1$;*

(iii) *$F_0 = B^3 \setminus F_1$ up to a negligible set.*

*Then, $\sup_{t\in[0,1]} P(F_t; B^3) \geq \pi + \delta_0$.*

*Proof.* Pick $\varepsilon = \mathscr{H}^3(B^3)/12 = \pi/9$ and consider $\delta_0 > 0$ to be the associated $\delta$ given by Lemma 3.6. If

$$\sup_{t\in[0,1]} P(F_t; B^3) < \pi + \delta_0,$$

then for every $t \in [0,1]$ there exists a $\mathbb{D}_n$-equivariant half-ball $\tilde{F}_t$ such that $\mathscr{H}^3(F_t \triangle \tilde{F}_t) \leq \pi/9$. Note that the $\mathbb{D}_n$-equivariant half-balls are the upper and the lower half-balls and, for $n = 2$, also the two half-balls bounded by the plane containing $\xi_0, \xi_1$ and the two half-balls bounded by the plane containing $\xi_0, \xi_2$. In any case we deduce that, for every $t \in [0,1]$, the $\mathbb{D}_n$-equivariant half-ball $\tilde{F}_t$ is uniquely determined. Therefore, by continuity of the family $\{F_t\}_{t\in[0,1]}$, $\tilde{F}_t$ must be constant, but this contradicts the assumption that $F_0$ is the complement of $F_1$ in $B^3$. □

**Corollary 3.8.** *Let $\{\Sigma_t\}_{t\in[0,1]}$ be the $\mathbb{D}_{g+1}$-sweepout given by Lemma 2.2 and let $\Pi$ be its $\mathbb{D}_{g+1}$-saturation given in Definition 3.2. Then the min-max width of $\Pi$ is larger than $\pi$ and smaller than $3\pi$, namely $\pi < W_\Pi < 3\pi$.*

## 4. Controlling the topology

In this section we prove Theorem 1.1. By Corollary 3.8, all conditions for applying Theorem 3.2 in [13], in particular the mountain-pass condition

$$W_\Pi > \max\{\mathscr{H}^2(\Sigma_0), \mathscr{H}^2(\Sigma_1)\},$$

are satisfied. We thus obtain a min-max sequence $\{\Sigma^j\}_{j\in\mathbb{N}}$ consisting of $\mathbb{D}_{g+1}$-equivariant surfaces and converging in the sense of varifolds to $m\Gamma$, where $\Gamma$ is a smooth, properly embedded, compact connected free boundary minimal surface in $B^3$, and the multiplicity $m$ is a positive integer. Moreover, the following statements hold.

(i) The surface $\Gamma$ contains the horizontal axes $\xi_1, \ldots, \xi_{g+1}$ and intersects $\xi_0$ orthogonally.

(ii) The integer $m$ is odd.

(iii) $m \operatorname{genus}(\Gamma) \leq g$.

*Remark* 4.1. Observe that statement (i) is a consequence of the $\mathbb{D}_{g+1}$-equivariance (cf. also Lemma 3.4 and Lemma 3.5 in [12]). Point (ii) follows from the invariance with respect to the rotation of angle $\pi$ around the axes $\xi_1, \ldots, \xi_{g+1}$ and its (self-contained) proof can be found at the end of Section 7.3 in [13]. Most importantly, for what concerns our application, we note that a weaker version of (iii) is sufficient, namely the inequality $\operatorname{genus}(\Gamma) \leq g$, which is given by Theorem 9.1 in [17] (based on [3, Theorem 0.6], written for the closed case).

**Lemma 4.2.** *The multiplicity $m$ is equal to $1$ and $\Gamma$ is not a (topological) disc.*

*Proof.* Fraser and Schoen (see [5, Theorem 5.4]) proved that any free boundary minimal surface in $B^3$ has area at least $\pi$. By varifold convergence of the min-max sequence, it holds $m\mathscr{H}^2(\Gamma) = W_\Pi$. By Corollary 3.8, $\pi < m\mathscr{H}^2(\Gamma) < 3\pi$ whence we conclude $m < 3$. In fact, $m = 1$ since $m$ must be odd by (ii). As a result, $\mathscr{H}^2(\Gamma) > \pi$, which implies that $\Gamma$ is not isometric to the equatorial disc. However, according to [19], the equatorial disc is the only possible free boundary minimal disc in $B^3$ up to ambient isometries. □

To control the boundary of $\Gamma$ we will first need the following elementary result, which applies to any proper equivariant surface in $B^3$, irrespective of minimality.

**Lemma 4.3.** *Let $\Gamma \subset B^3$ be any smooth, properly embedded, $\mathbb{D}_{g+1}$-equivariant surface which contains the horizontal axes $\xi_1, \ldots, \xi_{g+1}$. Then their endpoints $q_k := (\cos(k\pi/(g+1)), \sin(k\pi/(g+1)), 0)$ for $k \in \{1, \ldots, 2g+2\}$ are all contained in the same connected component of $\partial\Gamma$.*

*Proof.* Given $k \in \{1, \ldots, g+1\}$, there exists a connected component $\sigma$ of $\partial \Gamma$ containing the point $q_k$ because $\xi_k \subset \Gamma$ and $\Gamma$ is properly embedded. Let $\psi_k \in \mathbb{D}_{g+1}$ be the rotation of angle $\pi$ around $\xi_k$. Since $\partial \Gamma$ is $\mathbb{D}_{g+1}$-equivariant, we have in particular $\psi_k(\sigma) \subset \partial \Gamma$. In fact, $\psi_k(\sigma) = \sigma$ because $\sigma$ is a connected component of $\partial \Gamma$ intersecting $\psi_k(\sigma)$ at least in the point $q_k = \psi_k(q_k)$. Moreover, as $\Gamma$ is properly embedded and smooth, $\sigma \subset \partial B^3$ must be a smooth, simple closed curve. Any such curve divides $\partial B^3$ into two connected open domains $A'$ and $A''$. We then note that $\psi_k$ leaves the set $\partial A' = \sigma = \partial A''$ invariant, and that $A' = \psi_k(A'')$. It follows that $A'$ and $A''$ have the same area because $\psi_k$ is an isometry. Moreover, $A' = \psi_k(A'')$ implies that the antipodal point $q_{k+g+1} \in \xi_k$, which is fixed under $\psi_k$, must also be contained in $\sigma$ because it cannot be contained in $A'$ nor $A''$.

Now suppose, for the sake of a contradiction, that the point $q_\ell$ belongs to a different connected component $\varsigma$ of $\partial \Gamma$, for some $\ell \in \{1, \ldots, g+1\}$ with $\ell \neq k$. Then either $\varsigma \subset A'$ or $\varsigma \subset A''$ because $\sigma$ and $\varsigma$ are disjoint by definition. However, the whole argument in the previous paragraph applies to $\varsigma$ as well. Yet, in either case (i.e. both when $\varsigma \subset A'$ and $\varsigma \subset A''$), it is impossible that $\varsigma$ divides $\partial B^3$ into two domains of equal area, and this concludes the proof. □

**Lemma 4.4.** *The number of boundary components of $\Gamma$ is 1.*

*Proof.* Suppose, towards a contradiction, that $\partial \Gamma$ has more than one connected component. Then Lemma 4.3 implies that one connected component, say $\gamma$, of $\partial \Gamma$ is disjoint from $\mathcal{S} = \xi_0 \cup \xi_1 \cup \ldots \cup \xi_{g+1}$. Recall that the vertical axis $\xi_0$ is always disjoint from $\partial \Gamma$ because $\Gamma \subset B^3$ is properly embedded and intersects $\xi_0$ orthogonally by item (i) above. Moreover, let $\tilde{\gamma}$ be a simple closed curve in the interior of $\Gamma \setminus \mathcal{S}$ that is homotopic to $\gamma$ in $\Gamma \setminus \mathcal{S}$ (it is sufficient to slightly push $\gamma$ towards the interior of $\Gamma$).

Let then $\delta > 0$ be so small that

$$U_\delta \Gamma := \{x \in B^3 \,:\, \mathrm{dist}_{\mathbb{R}^3}(x, \Gamma) < \delta\}$$

is a tubular neighbourhood of $\Gamma$ in $B^3$. Since $\Gamma$ is connected and thus it is not a disc (since it has at least two boundary components), $\gamma$ and $\tilde{\gamma}$ are not contractible in $U_\delta \Gamma$.

Thanks to Proposition 4.10 in [17], without loss of generality we can assume that the min-max sequence $\{\Sigma^j\}_{j \in \mathbb{N}}$ is outer almost minimising (see [17, Definition 3.6]) in sufficiently small annuli. Therefore, by Simon's Lifting Lemma (see [3, Proposition 2.1] and [17, Section 9]), for every $j$ sufficiently large there exists a closed curve $\gamma^j \subset \Sigma^j \cap U_\delta \Gamma$ that is homotopic to $\tilde{\gamma}$ in $U_\delta \Gamma$ (note that we can apply the lemma since $\tilde{\gamma}$ is contained in the interior of $\Gamma$). In fact, given any $\rho > 0$, it follows from the proof in Section 4.3 of [3] that $\gamma^j$ can be taken in $U_\rho \tilde{\gamma} := \{x \in B^3 \,:\, \mathrm{dist}_{\mathbb{R}^3}(x, \tilde{\gamma}) < \rho\}$. Hence, choosing $\rho$ such that $U_\rho \tilde{\gamma} \subset U_{\delta/2} \Gamma \setminus \mathcal{S}$, we can guarantee that $\gamma^j \subset (\Sigma^j \cap U_{\delta/2} \Gamma) \setminus \mathcal{S}$ for every $j$ sufficiently large.

Now observe that $\Sigma^j$ is diffeomorphic to an element $\Sigma_{t_j}$ of the sweepout given by Lemma 2.2 for some $0 < t_j < 1$, through a $\mathbb{D}_{g+1}$-equivariant diffeomorphism. Thus, both connected components of $\Sigma^j \setminus (\xi_1 \cup \ldots \cup \xi_{g+1})$ are topological discs, hence $\gamma^j \subset \Sigma^j \setminus \mathcal{S}$ is contractible. Let us denote by $D^j$ a disc in the interior of $\Sigma^j$ such that $\partial D^j = \gamma^j$. We claim that $\gamma^j$ bounds a disc in $U_\delta \Gamma$ as well, which would contradict the existence of $\gamma$, since $\gamma^j$ is homotopic to $\tilde{\gamma}$ in $U_\delta \Gamma$ and $\tilde{\gamma}$ is not contractible there.

We will now exploit an argument similar to the one in Section 2.4 of [3]. Since the min-max sequence $\{\Sigma^j\}_{j \in \mathbb{N}}$ converges in the sense of varifolds to $\Gamma$, it follows that given any $\eta > 0$ there exists $J = J(\delta, \eta) \in \mathbb{N}$ such that, for every $j \geq J$,

$$\mathscr{H}^2(\Sigma^j \setminus U_{\delta/2} \Gamma) < \eta.$$

Defining $V_s\Gamma := \partial(U_s\Gamma) \cap B^3$ for $s \in {]}0,\delta[$, we observe that $\{V_s\Gamma\}_{s\in{]}0,\delta[}$ is a smooth foliation of $U_\delta\Gamma$ and we can apply the coarea formula to conclude

$$\int_{\delta/2}^{\delta} \mathscr{H}^1(\Sigma^j \cap V_s\Gamma) \leq \mathscr{H}^2(\Sigma^j \setminus U_{\delta/2}\Gamma) < \eta,$$

for every $j$ sufficiently large. Thus, there exists a subset $I \subset {]}\delta/2, \delta[$ of measure at least $\delta/4$ such that for all $s \in I$

$$\mathscr{H}^1(\Sigma^j \cap V_s\Gamma) < \frac{4\eta}{\delta}.$$

By Sard's theorem there exists $s \in I$ such that the intersection $\Xi_s^j := \Sigma^j \cap V_s\Gamma$ is transverse. This implies that any connected component of $\Xi_s^j$ is smooth and either a simple closed curve or an arc connecting two points of $\partial\Sigma^j$ in $V_s\Gamma$.

There exists $\lambda > 0$ (depending on $\Gamma$ and $\delta$) such that for any $s \in {]}\delta/2, \delta[$ any simple closed curve in $V_s\Gamma$ with length less than $\lambda$ bounds an embedded disc in $V_s\Gamma$. At this stage, we may choose $\eta > 0$ such that $4\eta < \lambda\delta$ and then $j \geq J(\delta, \eta)$ to ensure that the length of each connected component of $\Xi_s^j$ is less than $\lambda$. Now observe that $D^j \cap V_s\Gamma \subset \Sigma^j \cap V_s\Gamma = \Xi_s^j$. In particular, $D^j \cap V_s\Gamma$ consists of a finite number of simple closed curves (since $D^j$ is contained in the interior of $\Sigma^j$) of length less than $\lambda$ and thus each connected component of $D^j \cap V_s\Gamma$ bounds a disc in $V_s\Gamma$.

Hence, defining $G^j \subset U_\delta\Gamma$ as the connected component of $D^j \cap U_s\Gamma$ containing $\partial D^j = \gamma^j$, it is possible to cap the boundary components of $G^j$ lying in $V_s\Gamma$ with discs (contained in $U_\delta\Gamma$) in such a way that the resulting surface, which we denote by $\tilde{D}^j$, satisfies $\partial \tilde{D}^j = \partial D^j = \gamma^j$ (for recall that $\gamma^j \subset U_{\delta/2}\Gamma$). Note that $\tilde{D}^j$ is a topological disc since it is obtained from the topological disc $D^j$ by removing interior discs and then gluing discs with those same boundaries. Therefore it follows that $\tilde{\gamma}^j = \partial\tilde{D}^j$ is contractible in $U_\delta\Gamma$, which contradicts the initial choice of $\gamma$. $\square$

**Lemma 4.5.** *The genus of $\Gamma$ is $g$.*

*Proof.* By Lemmata 4.2 and 4.4, $\Gamma$ is not a topological disc and has connected boundary. Moreover, again by Lemma 4.2, $m = 1$ and, as a result (invoking fact (iii) in the weaker form needed for multiplicity one, cf. Remark 4.1), genus$(\Gamma) \leq g$. Because of this fact, and recalling that $\Gamma$ contains the origin due to item (i) above, Lemma B.1 applies, which proves the claim. $\square$

*Proof of Theorem 1.1.* According to Lemmata 4.4 and 4.5, the free boundary minimal surface $M_g = \Gamma$ has genus $g$ and connected boundary. As said before, $\Gamma$ inherits the dihedral symmetry $\mathbb{D}_{g+1}$ from the min-max sequence $\{\Sigma^j\}_{j\in\mathbb{N}}$, which completes the proof. $\square$

## A. Maximality of the unstable catenoid

In Section 2, we introduced the surfaces $C_s^{r,h}$, parametrised by $s \geq 0$, interpolating between the cylinder $C_0^{r,h}$ of radius $r$ and height $2h$ (for $s = 0$) and the union of two discs of radius $r$ with a line segment (as one lets $s \to \infty$). As we remarked, $C_s^{r,h}$ is minimal if $s$ is a solution to equation (8). Said $t_0 > 0$ the only positive solution of $\cosh(t) = t\sinh(t)$, equation (8) has two positive solutions $s_1$ and $s_2$ provided that $0 < h < r/\sinh(t_0)$ (note that the number $1/\sinh(t_0)$ is bounded from below by $0.6627$). In our notation, $s_2$ is the larger solution and corresponds to the so-called *unstable catenoid*. In this appendix we show that, as soon as we are willing to impose a slightly more restrictive condition on the ratio $h/r$, such a surface has largest area among all elements of the family in question. As a result, estimate (9), which plays an essential role in Section 2, follows at once from [14, Proposition 2.1].

**Lemma A.1.** *The unstable catenoid has largest area among all surfaces in the family $\{C_s^{r,h}\}_{s\geq 0}$ provided that $0 < 2h < r\tanh(1)$.*

*Proof.* As defined in (7), the surface $C_s^{r,h}$ is obtained by rotating the graph of $\rho\colon [-h, h] \to \mathbb{R}$ given by

$$\rho(z) = \frac{r\cosh(sz)}{\cosh(sh)} \tag{15}$$

around the vertical axis. Being a surface of revolution, the mean curvature of $C_s^{r,h}$ is easily computed to be

$$H = \frac{\rho\rho'' - (\rho')^2 - 1}{\left(1 + (\rho')^2\right)^{3/2}\rho}. \tag{16}$$

We notice that the denominator of (16) is strictly positive and that the numerator

$$\rho\rho'' - (\rho')^2 - 1 = \frac{r^2 s^2}{\cosh^2(sh)} - 1 \tag{17}$$

is independent of $z$. In particular, it follows that $H$ has the same sign of the function $rs - \cosh(sh)$. Recalling equation (8), the inequality $rs > \cosh(sh)$ is equivalent to $s_1 < s < s_2$ by strict convexity of $s \mapsto \cosh(sh)$. In this case $H > 0$, which implies for $s \in \,]s_1, s_2[$ the area of $C_s^{r,h}$ is an increasing function of $s$. Conversely, if $s > s_2$ or if $s < s_1$, then $H < 0$ and the area of $C_s^{r,h}$ is decreasing in $s$. This shows that the area of $C_s^{r,h}$ has a local minimum at $s_1$ and a local maximum at $s_2$. In order to prove the claim that $s_2$ is in fact a global maximum, it remains to check $\mathscr{H}^2(C_{s_2}^{r,h}) > \mathscr{H}^2(C_0^{r,h})$, i.e. that the unstable catenoid has larger area than the cylinder provided that $h/r$ is sufficiently small.

The area $A(s) = \mathscr{H}^2(C_s^{r,h})$ can be computed using the formula

$$A(s) = 4\pi \int_0^h \rho\sqrt{(\rho')^2 + 1}\, dz.$$

With $\rho$ as defined in (15), which in particular satisfies $\rho'' = s^2\rho$, the function $f = 4\pi\rho\sqrt{(\rho')^2 + 1}$ has a primitive given by

$$F = \frac{2\pi}{s^2}\left(\operatorname{asinh}(\rho') + \rho'\sqrt{(\rho')^2 + 1}\right).$$

Since $\rho'(0) = 0$ and $\rho'(h) = rs\tanh(sh)$, we obtain

$$A(s) = F(h) - F(0) = \frac{2\pi}{s^2}\left(\operatorname{asinh}(rs\tanh(sh)) + rs\tanh(sh)\sqrt{1 + r^2 s^2 \tanh^2(sh)}\right). \tag{18}$$

Comparing the derivatives of $s \mapsto rs$ and $s \mapsto \cosh(sh)$ at the intersection $s = s_2$, we obtain

$$r \leq h\sinh(s_2 h) < h\cosh(s_2 h) = rs_2 h$$

which implies $s_2 h > 1$. In turn, this yields

$$A(s_2) = \frac{2\pi}{s_2^2}\left(s_2 h + rs_2\tanh(s_2 h)\cosh(s_2 h)\right) = \frac{2\pi h}{s_2} + 2\pi r^2\tanh(s_2 h) > 2\pi r^2\tanh(1).$$

The assumption $2h < r\tanh(1)$ then implies $A(s_2) > 4\pi r h = \mathscr{H}^2(C_0^{r,h})$ which completes the proof. $\square$

# B. Structure of equivariant surfaces

**Lemma B.1.** *Let $1 \leq g \in \mathbb{N}$ and let $\Gamma \subset B^3$ be a compact, connected, properly embedded, $\mathbb{D}_{g+1}$-equivariant surface with genus $\gamma \in \{1, \ldots, g\}$ and one boundary component. Moreover, assume that $\Gamma$ contains the origin. Then $\Gamma$ has genus $\gamma = g$.*

*Proof.* In the case $g = 1$, there is nothing to prove. Therefore, let us assume $g \geq 2$. Since the given surface $\Gamma$ is $\mathbb{D}_{g+1}$-equivariant any intersection with the axis $\xi_0$ is orthogonal, hence (by embeddedness) $\Gamma \cap \xi_0$ shall consist of a finite set of points. Furthermore, since $\Gamma$ contains the origin, it will intersect the vertical axis $\xi_0$ in $2j+1$ many points, where $j$ is a nonnegative integer; also note that the poles $(0,0,1)$ and $(0,0,-1)$ cannot be contained in $\Gamma$ since it is properly embedded by assumption.

Let $C_{g+1} < \mathbb{D}_{g+1}$ be the cyclic subgroup of order $(g+1)$ which is generated by the rotation of angle $2\pi/(g+1)$ around $\xi_0$. The quotient $\Gamma' = \Gamma/C_{g+1}$ is a compact topological surface with boundary, i.e. a compact topological space in which every point has an open neighbourhood homeomorphic to some open subset of the (closed) half-plane. In particular, its Euler characteristic $\chi(\Gamma')$ is well-defined. Since the boundary of $\Gamma$ is connected and disjoint from the singular locus $\xi_0$ of the $C_{g+1}$-group action, the quotient $\Gamma'$ also has connected boundary homeomorphic to a circle. Moreover, $\Gamma'$ is orientable since every element of $C_{g+1}$ is orientation-preserving on $\Gamma$, which follows from the fact that the unit normal to $\Gamma$ at the origin is fixed under the action of $C_{g+1}$. Being orientable with connected boundary, $\Gamma'$ has Euler characteristic $\chi(\Gamma') = 1 - 2\gamma'$ for some integer $\gamma' \geq 0$. A suitable version of the Riemann–Hurwitz formula (see Remark B.2 below) implies

$$1 - 2\gamma = \chi(\Gamma) = (g+1)\chi(\Gamma') - (2j+1)g = (g+1)(1-2\gamma') - 2jg - g,$$

which is equivalent to $\gamma = (g+1)\gamma' + jg$. The assumption $\gamma \in \{1, \ldots, g\}$ then enforces $\gamma' = 0$ and $j = 1$ which proves $\gamma = g$. □

*Remark* B.2 (Riemann–Hurwitz formula, see e.g. Chapter IV.3 in [8]). Let $\Gamma$ and $\Gamma'$ be as in the proof of Lemma B.1. Let $T'$ be a triangulation of $\Gamma'$ such that every branch point of $\Gamma'$ is a vertex. Away from the branch points, the canonical projection $\Gamma \to \Gamma'$ is a covering map. Therefore, after refining $T'$ if necessary, the preimage of every triangle is a disjoint union of triangles in $\Gamma$ which leads to a triangulation $T$ of $\Gamma$. Note that $T$ has $(g+1)$ times as many faces and edges as $T'$. However, the $(2j+1)$ many branch points in $\Gamma'$ (which correspond to the points in $\Gamma \cap \xi_0$) have only one preimage rather than $(g+1)$. By the very definition of Euler characteristic we then have

$$\chi(\Gamma) = (g+1)\chi(\Gamma') - (2j+1)g.$$

It is appropriate to note that, since $\Gamma'$ is actually an orbifold, one could get to the same conclusions by invoking, in lieu of the Riemann–Hurwitz formula, a suitable version of the Gauss–Bonnet theorem for surfaces with conical singularities.

14 January 2020

Alessandro Carlotto
ETH D-Math, Rämistrasse 101, 8092 Zürich, Switzerland
*E-mail address:* `alessandro.carlotto@math.ethz.ch`

Giada Franz
ETH D-Math, Rämistrasse 101, 8092 Zürich, Switzerland
*E-mail address:* `giada.franz@math.ethz.ch`

Mario B. Schulz
Queen Mary University of London, School of Mathematical Sciences, Mile End Road, London E1 4NS, United Kingdom
*E-mail address:* `m.schulz@qmul.ac.uk`